\documentstyle[11pt,amscd,amssymb]{amsart}

\newcommand{\ar}{\rightarrow}
\newcommand{\bd}{\partial}
\newcommand{\x}{\times}
\newcommand{\ox}{\otimes}
\newcommand{\iso}{\cong}
\newcommand{\isom}{\stackrel{\simeq}{\ar}}

\newcommand{\la}{\langle}
\newcommand{\ra}{\rangle}

\newcommand{\RR}{{\Bbb R}}
\newcommand{\TT}{{\Bbb T}}
\newcommand{\SS}{{\Bbb S}}
\newcommand{\ZZ}{{\Bbb Z}}

\renewcommand{\a}{\alpha}
\renewcommand{\b}{\beta}
\renewcommand{\d}{\delta}
\newcommand{\g}{\gamma}

\renewcommand{\l}{\lambda}

\newcommand{\G}{\Gamma}

\theoremstyle{plain}
\begingroup
\newtheorem{thm}{Theorem}

\newtheorem{lem}[thm]{Lemma}

\endgroup

\theoremstyle{definition}

\theoremstyle{remark}

\title{On non-formal simply connected manifolds}

\author{Marisa Fern\'{a}ndez}
\address{ Departamento de Matem\'aticas \\
Universidad del Pa\'{\i}s Vasco \\ Apartado 644 \\ 48080 Bilbao \\
Spain} \email{mtpferol@@lg.ehu.es}

\author{Vicente Mu\~noz}
\address{Departamento de Matem\'aticas \\
Universidad Aut\'onoma de Madrid
\\ 28049 Madrid \\ Spain}
\email{vicente.munoz@@uam.es}

\thanks{\noindent First author partially supported by CICYT
(Spain) Project BFM2001-3778-C03-02 and UPV
00127.310-E-14813/2002. \\ Second author supported by CICYT
project BFM2000-0024. \\ Also partially supported by The European
Contract Human Potential Programme, Research Training Network
HPRN-CT-2000-00101.}

\date{January, 2003.}

\begin{document}

\begin{abstract}
We construct examples of non-formal simply connected and compact
oriented manifolds of any dimension bigger or equal to $7$.
\end{abstract}

\maketitle

\section{Introduction}

An oriented compact manifold of dimension at most $2$ is formal.
On the other hand, if the dimension is $3$ or more, there are
examples which are non-formal, e.g.,\ nilmanifolds which are not
tori \cite{Hasegawa}.

If we turn our attention to simply connected manifolds, we know
that a simply connected oriented compact manifold of dimension at
most $6$ is formal \cite{NM, Mi, FM}. The natural question already
raised in \cite{FM} is whether there are examples of non-formal
simply connected oriented compact manifolds of dimension $d\geq
7$.

Clearly, the question is reduced to the cases $d=7$ and $d=8$. For
if we have a non-formal simply connected manifold $M$ of dimension
$d$, then $M\times S^{2n}$ is a non-formal simply connected
manifold of dimension $d+2n$, for any $n \geq 1$.

{}From now on let $d=7$ or $d=8$. By the results of \cite{FM}, if
a $d$--dimensional connected and compact oriented manifold $M$ is
$3$--formal then it is formal. Therefore, the non-formality of $M$
has to be detected in the $3$--stage of its minimal model.
Moreover if $H^1(M)=0$ then $M$ is automatically $2$--formal, so
the non-formality is due to the kernel of the cup product map
$\cup: H^2(M) \ox H^2(M) \to H^4(M)$. The easiest way to detect
the non-formality is thus to have a non-trivial Massey product of
cohomology classes of degree $2$.

The method of construction of $d$--dimensional simply connected
manifolds that we will use is the following: take a non-formal
compact nilmanifold $X$ of dimension $d$ with a non-trivial Massey
product of cohomology classes of degree $1$. Multiply these
cohomology classes by some cohomology classes so that we get a
non-trivial Massey product of cohomology classes of degree $2$.
Then perform a suitable surgery of $X$ to kill the fundamental
group such that the non-trivial Massey product survives. This will
give the sought example.

In \cite{BT} Babenko and Taimanov have already given examples of
non-formal simply connected manifolds of any {\em even\/}
dimension bigger or equal to $10$. The relevant property of their
examples is that they are symplectic manifolds. They ask whether
there exist examples of non-formal simply connected {\em
symplectic\/} manifolds of dimension $8$. Unfortunately, our
examples do not have a symplectic structure, at least in an
obvious way.

\section{The $8$-dimensional example}

Let $H$ be the Heisenberg group, that is, the connected nilpotent
Lie group of dimension $3$ consisting of matrices of the form
 $$
 a=\left( \begin{array}{ccc} 1&x&z\\ 0&1&y\\ 0&0&1 \end{array}\right),
 $$
where $x,y,z \in {\RR}$. Then a global system of coordinates
${x,y,z}$ for $H$ is given by $x(a)=x$, $y(a)=y$, $z(a)=z$, and a
standard calculation shows that a basis for the left invariant
$1$--forms on $H$ consists of $\{dx, dy, dz-xdy\}$. Let $\Gamma$
be the discrete subgroup of $H$ consisting of matrices whose
entries are integer numbers. So the quotient space $N
=\Gamma{\backslash}H$ is a compact $3$--dimensional nilmanifold.
Hence the forms $dx$, $dy$, $dz-xdy$ descend to $1$--forms
$\alpha$, $\beta$, $\gamma$ on $N$ and
 $$
 d\alpha=d\beta=0, \quad  d\gamma=-\alpha \wedge \beta.
 $$
The non-formality of $N$ is detected by a non-zero triple Massey
product
 $$
 \la [\a],[\b],[\a]\ra= [2 \, \a \wedge \g].
 $$

Now let us consider $X=N \x \TT^5$, where $\TT^5=\RR^5/\ZZ^5$. The
coordinates of $\RR^5$ will be denoted $x_1,x_2,x_3,x_4,x_5$. So
$\{dx_i| 1\leq i \leq 5\}$ defines a basis $\{\d_i |1\leq i \leq
5\}$ for the $1$--forms on $\TT^5$. By multiplying the classes
$\a$ and $\b$ by some of the $\d_i$, we get a non-zero triple
Massey product of cohomology classes of degree $2$ for $X$,
 \begin{equation} \label{eqn:massey}
 \la [\a\wedge \d_1],[\b\wedge \d_2],[\a\wedge \d_3]\ra=
 [2\, \g \wedge \a\wedge \d_1\wedge \d_2\wedge \d_3].
 \end{equation}

Our aim now is to kill the fundamental group of $X$ by performing
a suitable surgery construction. Let  $C_1$ the image of
$\{(x,0,0)|x\in\RR\}\subset H$ in $N=\G\backslash N$ and let $C_2$
be the image of $\{(0,y,\xi)| y\in\RR\}$ in $N$, where $\xi$ is a
generic real number. Then $C_1,C_2\subset N$ are disjoint embedded
circles such that $p(C_1)=\SS^1 \x \{0\}$, $p(C_2)=\{0\}\x\SS^1$.
The projection $p(x,y,z)=(x,y)$ describes $N$ as a fiber bundle
$p:N\to \TT^2$ with fiber $\SS^1$. Actually, $N$ is the total
space of the unit circle bundle of the line bundle of degree $1$
over the $2$--torus. The fundamental group of $N$ is therefore
  \begin{equation}\label{eqn:pi}
  \pi_1(N) \iso \Gamma = \la \l_1,\l_2,\l_3 | [\l_1,\l_2]=\l_3,
    \text { $\l_3$ central} \ra,
  \end{equation}
where $\l_3$ corresponds to the fiber, $\l_1$ and $\l_2$
correspond to the homotopy classes $\l_1=[C_1]$ and $\l_2=[C_2]$.
The fundamental group of $X=N\x\TT^5$ is
 \begin{equation}\label{eqn:pi1}
 \pi_1(X)=\pi_1(N) \oplus \ZZ^5.
 \end{equation}

Consider the following submanifolds embedded in $X$:
\begin{eqnarray*}
 T_1&=& C_1 \x \SS^1 \x \{0\} \x \SS^1 \x \{0\} \x \SS^1, \\
 T_2&=& C_2 \x \{0\} \x \SS^1 \x \{0\} \x \SS^1 \x \SS^1,
\end{eqnarray*}
which are $4$-dimensional tori with trivial normal bundle.
Consider now another $8$-manifold $Y$ with an embedded
$4$-dimensional torus $T$ with trivial normal bundle. Then we may
perform the {\em fiber connected sum\/} of $X$ and $Y$ identifying
$T_1$ and $T$, denoted $X \#_{T_1=T} Y$, in the following way:
take (open) tubular neighborhoods $\nu_1 \subset X$ and $\nu
\subset Y$ of $T_1$ and $T$ respectively; then $\bd \nu_1 \iso
{\TT^4}\x \SS^3$ and $\bd \nu \iso {\TT^4}\x \SS^3$; take an
orientation reversing diffeomorphism $\phi:\bd \nu_1 \isom \bd
\nu$; the fiber connected sum is defined to be the (oriented)
manifold obtained by gluing $X-\nu_1$ and $Y-\nu$ along their
boundaries by the diffeomorphism $\phi$. In general, the resulting
manifold depends on the identification $\phi$, but this will not
be relevant for our purposes.

\begin{lem}\label{lem:pi1}
  Suppose $Y$ is simply connected. Then the fundamental group of
  $X \#_{T_1=T} Y$
  is the quotient of $\pi_1(X)$ by the image of $\pi_1(T_1)$.
\end{lem}

\begin{pf}
 Since the codimension of $T_1$ is bigger or equal than $3$, we
 have that
 $\pi_1(X-\nu_1)=\pi_1(X-T_1)$ is isomorphic to $\pi_1(X)$. The
 Seifert-Van Kampen theorem establishes that $\pi_1(X \#_{T_1=T} Y)$ is
 the amalgamated sum of $\pi_1(X-\nu_1)=\pi_1(X)$ and
 $\pi_1(Y-\nu)=\pi_1(Y)=1$ over the image of $\pi_1(\bd
 \nu_1)=\pi_1(T_1\x\SS^3)=\pi_1(T_1)$, as required.
\end{pf}

We shall take for $Y$ the sphere $\SS^8$. We embed a
$4$-dimensional torus $\TT^4$ in $\RR^8$. This torus has a trivial
normal bundle since its tangent bundle is trivial (being
parallelizable) and the tangent bundle of $\RR^8$ is also trivial.
After compatifying $\RR^8$ by one point we get a $4$-dimensional
torus $T\subset \SS^8$ with trivial normal bundle.

In the same way, we may consider another copy of the
$4$-dimensional torus $T\subset \SS^8$ and perform the fiber
connected sum of $X$ and $\SS^8$ identifying $T_2$ and $T$. We may
do both fiber connected sums along $T_1$ and $T_2$ simultaneously,
since $T_1$ and $T_2$ are disjoint. Call
  $$
  M= X \#_{T_1=T} \SS^8 \#_{T_2=T} \SS^8
  $$
the resulting manifold. By Lemma \ref{lem:pi1}, $\pi_1(M)$ is the
quotient of $\pi_1(X)$ by the images of $\pi_1(T_1)$ and
$\pi_1(T_2)$. This kills the $\ZZ^5$ summand in \eqref{eqn:pi1}
and it also kills $\l_1$ and $\l_2$ in \eqref{eqn:pi}. Therefore
$\pi_1(M)=1$, i.e.,\, $M$ is simply connected.

\section{Non-formality of the constructed manifold}

Our goal is now to prove that $M$ is non-formal. We shall do this
by proving the non-vanishing of a suitable triple Massey product.
More specifically, let us prove that the Massey product
\eqref{eqn:massey} survives to $M$. For this, let us describe
geometrically the cohomology classes $[\a\wedge \d_1]$, $[\b\wedge
\d_2]$ and $[\a\wedge \d_3]$. Consider the following three
codimension $2$ submanifolds of $X$:
\begin{eqnarray*}
 B_1&=&p^{-1}(\SS^1\x\{a_1\})\x\{b_1\}\x \SS^1 \x \SS^1 \x\SS^1\x\SS^1,\\
 B_2&=&p^{-1}(\{a_2\}\x\SS^1)\x \SS^1 \x\{b_2\}\x \SS^1 \x\SS^1\x\SS^1,\\
 B_3&=&p^{-1}(\SS^1\x\{a_3\})\x \SS^1 \x \SS^1 \x\{b_3\}\x\SS^1\x\SS^1,
\end{eqnarray*}
where the $a_i$ and $b_i$ are generic points of $\SS^1$. It is
easy to check that $B_i\cap T_j=\emptyset$ for all $i$ and $j$. So
$B_i$ may be also considered as submanifolds of $M$. Let $\eta_i$
be the $2$--forms representing the Poincar\'{e} dual to $B_i$ in $X$.
By \cite{BoTu}, $\eta_i$ are taken supported in a small tubular
neighborhood of $B_i$. Therefore the support of $B_i$ lies inside
$X-T_1-T_2$, so we also have naturally $\eta_i \in \Omega^2(M)$.
Note that in $X$ we have clearly that $[\eta_1]=[\a \wedge e_1]$,
$[\eta_2]= [\b \wedge e_2]$ and $[\eta_3]=[\a\wedge e_3]$, where
$e_i$ are differential $1$--forms on $\SS^1$ cohomologous to
$\d_i$ and supported in a neighborhood of $b_i\in \SS^1$. Thus
$[\eta_1]=[\a \wedge \d_1]$, $[\eta_2]= [\b \wedge \d_2]$ and
$[\eta_3]=[\a\wedge \d_3]$.

\begin{lem} \label{Massey}
 The triple Massey product $\la [\eta_1],[\eta_2],[\eta_3]\ra$ is
 well-defined on $M$ and equals to $[2\, \g\wedge\a\wedge e_1\wedge
 e_2 \wedge e_3]$.
\end{lem}

\begin{pf}
Clearly
 $$
 (\a \wedge e_1) \wedge (\b \wedge e_2) = d\g \wedge e_1\wedge
 e_2,
 $$
where the $3$--form $\g\wedge e_1\wedge e_2$ is supported in a
neighborhood of $N \x  \{b_1\} \x\{b_2\} \x \SS^1 \x\SS^1 \x
\SS^1$, which is disjoint from $T_1$ and $T_2$. Hence $\g \wedge
e_1\wedge e_2$ is well-defined as a form in $M$. Also
 $$
 (\b \wedge e_2) \wedge (\a \wedge e_3) = -d\g \wedge e_2\wedge
 e_3,
 $$
where $-\g \wedge e_2\wedge e_3$ is also well-defined in $M$. So
the triple Massey product
 $$
 \la [\eta_1],[\eta_2],[\eta_3]\ra =[2\,\g \wedge \a \wedge e_1\wedge
 e_2 \wedge e_3]
 $$
is well-defined in $M$.
\end{pf}

Finally let us see that this Massey product $\la
[\eta_1],[\eta_2],[\eta_3]\ra =[2\,\g \wedge \a \wedge e_1\wedge
e_2 \wedge e_3]$ is non-zero in
 $$
 \frac{H^5(M)}{[\a\wedge e_1] \cup H^3(M) + H^3(M) \cup [\a \wedge
 e_3]}.
 $$
To see this, consider $B_4=p^{-1}(\{a_4\} \x \SS^1) \x \SS^1 \x
\SS^1 \x \SS^1 \x \{b_4\} \x  \{b_5\}$, for generic points
$a_4,b_4,b_5$ of $\SS^1$. Then the Poincar\'{e} dual of $B_4$ is
defined by a $3$--form $\beta'\wedge e_4\wedge e_5$ supported near
$B_4$, where $\beta'$ is Poincar\'{e} dual to $p^{-1}( \{a_4\} \x
\SS^1)$ and $[\beta']=[\beta]$, $[e_4]=[\d_4]$ and $[e_5]=[\d_5]$.
Again this $3$--form can be considered as a form in $M$. Now for
any $[\varphi], [\varphi']\in H^3(M)$ we have
 $$
  ([2\, \g \wedge \a \wedge e_1\wedge e_2
 \wedge e_3]+ [\a\wedge e_1 \wedge \varphi] +[\b \wedge
 e_3\wedge \varphi']) \cdot [\beta'\wedge e_4 \wedge e_5]=-2,
 $$
since the first product gives $2$; to compute the second product,
we notice that the $5$--form $\a\wedge \b' \wedge e_1\wedge
e_4\wedge e_5$ is exact in $M$ because $\a\wedge \b' \wedge
e_1\wedge e_4\wedge e_5=-d\gamma' \wedge e_1\wedge e_4\wedge e_5$
in $X$, with $\gamma'=\gamma+f\,\alpha$ for some function $f$ on
$N$, and $\gamma'\wedge e_1\wedge e_4\wedge e_5$ is well-defined
on $M$; and for the third product, $\a\wedge \b' \wedge e_3\wedge
e_4\wedge e_5$ is also exact in $M$. Therefore we have proved the
following

\begin{thm}
 $M$ is a compact oriented simply connected non-formal $8$--manifold.
\end{thm}

\section{The $7$-dimensional example}

A compact oriented simply connected non-formal manifold $M'$ of
dimension $7$ is obtained in an analogous fashion to the
construction of the $8$--dimensional manifold $M$. We start with
$X'=N\x \TT^4$ and consider the $3$-dimensional tori
\begin{eqnarray*}
 T_1'&=& C_1 \x \SS^1 \x \{0\} \x \SS^1 \x \{0\}, \\
 T_2'&=& C_2 \x \{0\} \x \SS^1 \x \{0\} \x \SS^1.
\end{eqnarray*}
Define
 $$
 M'=X' \#_{T_1'=T'} \SS^7 \#_{T_2'=T'}\SS^7
 $$
where $T'$ is an embedded $3$--torus in $\SS^7$ with trivial
normal bundle. Then $M'$ is a non-formal simply connected
manifold. To prove the non-formality, consider the codimension $2$
submanifolds
\begin{eqnarray*}
 B_1'&=& p^{-1}(\SS^1 \x \{a_1\}) \x\{b_1\}\x \SS^1 \x \SS^1 \x \SS^1\\
 B_2'&=& p^{-1}(\{a_2\} \x \SS^1) \x \SS^1 \x\{b_2\}\x \SS^1 \x \SS^1\\
 B_3'&=& p^{-1}(\SS^1 \x \{a_3\}) \x \SS^1 \x \SS^1 \x\{b_3\}\x \SS^1
\end{eqnarray*}
and the $2$--forms $\eta_i'$ Poincar\'{e} dual to $B_i$. Then $\la
[\eta'_1],[\eta'_2],[\eta'_3]\ra =[2\,\g \wedge \a \wedge
e_1\wedge e_2 \wedge e_3]$. This triple Massey product is non-zero
in
 $$
 \frac{H^5(M')}{[\a\wedge e_1] \cup H^3(M') + H^3(M') \cup [\a \wedge
 e_3]},
 $$
by using the same argument as before with $B_4'=p^{-1}(\{a_4\} \x
\SS^1) \x \SS^1 \x \SS^1 \x \SS^1 \x \{b_4\}$.

Note that it is in this last step where the similar argument for
the $6$--dimensional case breaks down, since if we drop the last
factor all throughout the argument, then the submanifold
$B_4''=p^{-1}(\{a_4\} \x \SS^1) \x \SS^1 \x \SS^1 \x \SS^1$ would
not be disjoint from the two tori where the surgery is taken
place.


\begin{thebibliography}{10}

\bibitem{BT}
  I.K. Babenko, I.A. Taimanov, On nonformal simply connected
     symplectic manifolds,
     {\it Siberian Math. J.\/} {\bf 41} (2000), 204--217.

\bibitem{BoTu} R. Bott, L.W. Tu, {\em Differential forms in
  algebraic topology,\/} Graduate Texts in Maths, Vol.\ 82,
  Springer-Verlag, 1982.

\bibitem{FM}
  M. Fern\'{a}ndez, V. Mu\~{n}oz, On the formality and hard
  Lefschetz property for Donaldson symplectic manifolds,
     Preprint {\tt math.SG/0211017}.

\bibitem{Hasegawa}
 K. Hasegawa, Minimal models of nilmanifolds, {\it Proc. Amer.
Math. Soc.\/} {\bf 106} (1989), 65--71.


\bibitem{Mi}
 T.J. Miller, On the formality of $(k-1)$ connected compact manifolds
   of dimension less than or equal to $(4k-2)$,
   {\it Illinois. J. Math.\/} {\bf 23} (1979), 253--258.

\bibitem{NM}
 J. Neisendorfer, T.J. Miller, Formal and coformal spaces,
   {\it Illinois. J. Math.\/} {\bf 22} (1978), 565--580.


\end{thebibliography}
\end{document}